\documentclass[11pt]{article}
\usepackage{amssymb, amsmath, url, graphicx, setspace, geometry}
\usepackage{calrsfs}
\usepackage{wasysym}
\usepackage{xcolor}

\def\3{\subset }
\def\4{\subseteq }
\def\<{\left<}
\def\>{\right>}

\def\bit{\begin{itemize}}
\def\eit{\end{itemize}}
\def\3{\subset }
\def\4{\subseteq }

\def\0{\leqno}

\def\barr{\begin{array}}
\def\earr{\end{array}}

\def\Z{{\rlap{$\kern2pt{\rm Z}$}{\rm Z}\,}}

\title{\bf Groups whose same-order types are arithmetic progressions}
\author{Mihai-Silviu Lazorec and Marius T\u arn\u auceanu}
\date{}

\begin{document}
\maketitle

\begin{abstract}
The same-order type $\tau_e(G)$ of a finite group $G$ is a set formed of the sizes of the equivalence classes containing the same order elements of $G$. In this paper, we study an arithmetical property of this set. More exactly, we outline some results on the classification and existence of finite groups whose same-order types are arithmetic progressions formed of 3 or 4 elements, the latter being the maximum size of such a sequence. 
\end{abstract}



\section{Introduction}

Let $G$ be a finite group and let $n\geq 2$ be an integer. We denote the order of an element $x$ of $G$ and the cyclic group of order $n$ by $o(x)$ and $C_n$, respectively. The study of the element orders of $G$ is a very popular topic nowadays. During the last years, a lot of results concerning the nature/structure of $G$ were described via different tools involving its element orders. 
In this paper, we focus on the same-order type of $G$, i.e. the set of sizes of the equivalence classes with respect to an equivalence relation $``\sim"$ defined as follows:
$$\forall \ x, y\in G, x\sim y \Longleftrightarrow o(x)=o(y).$$

We denote the same-order type of $G$ by $\tau_e(G)$, as it was suggested in \cite{12}. Obviously, $\tau_e(G)=\lbrace 1\rbrace$ if and only if $G$ is trivial or $G\cong C_2$. Small sizes of same-order types may be used to detect the nilpotency or the solvability of finite groups. More exactly, in \cite{11} (see Theorems 2.1 and 3.7), Shen proved the following two main results:\\

\textbf{Theorem 1.1.} \textit{Let $G$ be a finite group such that $\tau_e(G)=\lbrace 1, n\rbrace$. Then $G$ is nilpotent and $G$ is isomorphic to one of the following groups:
\begin{itemize}
\item[(1)] a $p$-group $P$ with $|P|\geq 3$ and $exp(P)=p$;
\item[(2)] the quaternion group $Q_8$;
\item[(3)] the cyclic group $C_4$;
\item[(4)] a direct product $C_2\times P_1$, where $P_1$ is a $p$-group with $exp(P_1)=p$, for an odd prime $p$.
\end{itemize}}

\textbf{Theorem 1.2.} \textit{Let $G$ be a finite group such that $\tau_e(G)=\lbrace 1, m, n\rbrace$. Then $G$ is solvable and $G$ is isomorphic to one of the following groups:
\begin{itemize}
\item[(1)] a $p$-group with same-order type $\lbrace 1, m, n\rbrace$;
\item[(2)] the cyclic group $C_{12}$;
\item[(3)] the special linear group $SL(2, 3)$;
\item[(4)] the direct product $C_7\times Q_8$;
\item[(5)] a semidirect product $C_p^t\rtimes C_4$, where $p$ is an odd prime and $t\geq 1$ is an integer;
\item[(6)] a Frobenius group $C_2^t\rtimes C_p$, where $p$ is an odd prime and $t\geq 2$ is an integer;
\item[(7)] a Frobenius group $C_p^t\rtimes C_2$, where $p$ is an odd prime and $t\geq 1$ is an integer;
\item[(8)] a Frobenius group $P\rtimes C_q$, where $P$ is a $p$-group with $exp(P)=p$ and $p, q$ are distinct odd primes;
\item[(9)] a direct product $C_2\times H$, where $H$ is of type (1), (7) or (8); if $H$ is of type (1), $p$ is an odd prime.
\end{itemize}}

It is known that the order $|G|$ and the set $\pi_e(G)$ formed of the element orders of a finite group $G$ can be used together in order to recognize finite simple groups (Question 12.39 in the Kourovka Notebook \cite{10}). The same-order type of a finite group $G$ may be also connected with finite non-abelian simple groups. For instance, the main result of \cite{12} states that $\tau_e(G)=\lbrace 1, 15, 20, 24\rbrace$ if and only if $G\cong A_5$. Also, in \cite{9}, Jafari Taghvasani and Zarrin proved that $A_5$ is the only non-abelian simple group having a same-order type formed of 4 elements. 

Our aim is to study an arithmetical property of the same-order type of a finite group. More exactly, we are interested in providing some results concerning the classification and existence of finite groups $G$ for which $\tau_e(G)$ is an arithmetic progression. We show that the maximum size of such a same-order type is 4. Our main results concern the same-order types formed of 3 or 4 elements. More exactly, we prove that if $G$ is not a finite 2-group having more than one cyclic subgroup and $|\tau_e(G)|=3$, then $\tau_e(G)$ is an arithmetic progression if and only if $G\cong S_3$. Also, we show that there are no finite nilpotent groups whose same-order types are arithmetic progressions formed of 4 elements. 

We end this section by mentioning that several older and recent papers concerning arithmetical properties of some finite sets formed of different quantities related to finite groups were written. For instance, in \cite{3, 4}, Bianchi et al. classified all finite groups whose conjugacy class sizes form an arithmetic progression. Also, let us denote by $\pi_s(G)$ the set of proper subgroup orders and by $\pi_{as}(G)$ the set of proper abelian subgroup orders of a finite group $G$. We recall that the main results of Brandl and Shi \cite{5}, Shi \cite{13}, Feng \cite{6} and T\u arn\u auceanu \cite{14} highlighted all finite groups $G$ such that $\pi_e(G)$, $\pi_s(G)$ and $\pi_{as}(G)$ are arithmetic progressions.

\section{Main results}

Let $G$ be a finite group and let $n$ be a positive integer such that $n||G|$. We denote by $s_n$ and $c_n$ the number of elements of order $n$ and the number of cyclic subgroups of order $n$ of $G$, respectively. Also, Euler's totient function is denoted by $\varphi$. We recall some preliminary results through the following lemma.\\

\textbf{Lemma 2.1.} \textit{Let $G$ be a finite group. Then the following statements hold:
\begin{itemize}
\item[(i)] (see Section 2 of \cite{7}) If $n$ is a positive integer such that $n||G|$, then $n|\sum\limits_{m|n}s_m$.
\item[(ii)] (see Theorem 1 of \cite{8}) If $p$ is a prime number such that $p||G|$, then $s_p\equiv -1 \ (mod \ p)$.
\item[(iii)] (see Lemma 3.1 of \cite{11}) If $|G|$ is odd, $|\tau_e(G)|=3$ and $p, q$ are distinct prime divisors of $|G|$, then $s_p\ne s_q$.
\item[(iv)] (see Lemma 3 of \cite{12}) If $|G|\geq 2$ and $s$ is the maximal number of same order elements in $G$, then $|G|\leq s(s^2-1).$
\end{itemize}}

If one wants to find a finite group $G$ whose same-order type $\tau_e(G)$ is an arithmetic progression, the first step is to determine the maximum size of $\tau_e(G)$. Hence, the following proposition is a key result for our study since it answers to the previous question.\\

\textbf{Proposition 2.2.} \textit{Let $G$ be a finite group such that its same-order type $\tau_e(G)$ is an arithmetic progression. Then $|\tau_e(G)|\leq 4$.}\\

\textbf{Proof.} Let $G$ be a finite group such that $\tau_e(G)$ is an arithmetic progression of ratio $r\geq 1$. Suppose that $|\tau_e(G)|\geq 5$. Then $\tau_e(G)$ contains the odd integers $1+2r$ and $1+4r$. Hence, there are two distinct divisors, say $m>1$ and $n>1$, of $|G|$ such that $s_m=1+2r$ and $s_n=1+4r$. Moreover, without losing generality, we may suppose that $n\geq 3$. In this case $\varphi(n)$ is an even integer and we know that $\varphi(n)|s_n$. Therefore, $s_n$ would be an even integer, a contradiction. Consequently, $|\tau_e(G)|\leq 4$.
\hfill\rule{1,5mm}{1,5mm}\\

As a consequence of our last result, we can say that if $G$ is trivial or $G\cong C_2$, then $\tau_e(G)=\lbrace 1\rbrace$ is an arithmetic progression if we take into consideration the possibility of setting a null ratio. Also, it is obvious that all finite groups listed in Theorem 1.1 have the same arithmetical property for different positive ratios. More exactly, we have the same-order types $\tau_e(P)=\lbrace 1, |P|-1\rbrace$, $\tau_e(Q_8)=\lbrace 1, 6\rbrace$, $\tau_e(C_4)=\lbrace 1, 2\rbrace$ and $\tau_e(C_2\times P_1)=\lbrace 1, |P_1|-1\rbrace$. According to Proposition 2.2, it remains to see if there are any finite groups $G$ whose same-order types $\tau_e(G)$ are arithmetic progressions containing 3 or 4 elements. In what follows, we study these two remaining scenarios.\\

\textbf{Theorem 2.3.} \textit{Let $G$ be a finite group such that $G$ is not a 2-group with $c_2>1$ and $|\tau_e(G)|=3$. Then $\tau_e(G)$ is an arithmetic progression if and only if $G\cong S_3$}.\\

\textbf{Proof.} Let $G$ be a finite group such that $|\tau_e(G)|=3$. Then $G$ is isomorphic to one of the groups listed in Theorem 1.2. Moreover, we suppose that $\tau_e(G)$ is an arithmetic progression of ratio $r\geq 1$, so $\tau_e(G)=\lbrace 1, 1+r, 1+2r\rbrace$. We need to analyse 9 cases corresponding to the classification given by Theorem 1.2. 

\textit{Case 1.} Assume that $G\cong P$, where $P$ is chosen such that it is not a 2-group having more than one cyclic subgroup of order 2. If $P$ is a 2-group with an unique cyclic subgroup of order 2, then, according to Proposition 1.3 of \cite{2}, we have $G\cong C_{2^k}$ or $G\cong Q_{2^k}$, where $k\geq 3$ is an integer and $$Q_{2^k}=\langle x, y \ | \ x^{2^{k-1}}=y^4=1, x^y=x^{2^{k-1}-1}\rangle$$ is the generalized quaternion group. If $G\cong C_{2^k}$, then $\tau_e(G)=\lbrace 1, s_4=\varphi(4), \ldots, s_{2^k}=\varphi(2^k)\rbrace$. Since $|\tau_e(G)|=3$, we have $G\cong C_8$, so $\tau_e(G)=\lbrace 1, 2, 4\rbrace$ which is not an arithmetic progression, a contradiction. If $G\cong Q_{2^k}$, we may suppose that $k\geq 4$ since $|\tau_e(Q_8)|=2$. We have $\tau_e(G)=\lbrace 1, s_4=2(1+2^{k-2}), s_8=\varphi(8),\ldots,s_{2^{k-1}}=\varphi(2^{k-1})\rbrace$. Our hypotheses lead to $G\cong Q_{16}$ and $\tau_e(G)=\lbrace 1, 4, 10\rbrace$, a contradiction. 

If $P$ is a $p$-group and $p$ is odd, we may assume that $exp(G)\geq p^2$ since $|\tau_e(G)|=3$. We have $s_p=\varphi(p)c_p$ and $s_{p^2}=\varphi(p^2)c_{p^2}$. Then $p|s_{p^2}$ and, since $s_p\equiv -1 \ (mod \ p)$ by Lemma 2.1 (ii), the same-order type of $G$ is $\tau_e(G)=\lbrace 1, s_p, s_{p^2}\rbrace$. But both $s_p$ and $s_{p^2}$ are even. Therefore, $\tau_e(G)$ is not an arithmetic progression, a contradiction.

\textit{Cases 2, 3 and 4.} Assume that $G\cong C_{12}$, $G\cong SL(2,3)$ or $G\cong C_7\times Q_8$. Then $\tau_e(G)=\lbrace 1, 2, 4\rbrace, \tau_e(G)=\lbrace 1, 6, 8\rbrace$ or $\tau_e(G)=\lbrace 1, 6, 36\rbrace$, respectively. None of these is an arithmetic progression, a contradiction.

\textit{Case 5.} Suppose that $G\cong C_p^t\rtimes C_4$, where $p$ is an odd prime and $t\geq 1$. Then $s_p=p^t-1$ and $s_4=\varphi(4)c_4$. Since $c_4>1$, according to Sylow's 3rd Theorem, we obtain $p|c_4$, so $p|s_4$. Then $s_p\ne s_4$ and the same-order type of $G$ is $\tau_e(G)=\lbrace 1, s_4, s_p\rbrace$. However, both $s_4$ and $s_p$ are even. Consequently, $\tau_e(G)$ is not an arithmetic progression, a contradiction.

\textit{Case 6.} Let $p$ be an odd prime and $t\geq 2$. Suppose that $G\cong C_2^t\rtimes C_p$. Since $s_2=2^t-1$ is odd, $s_p=\varphi(p)c_p$ is even and $\tau_e(G)=\lbrace 1, s_2, s_p\rbrace$ is an arithmetic progression, we deduce that $s_p=1+r$ and $s_2=1+2r$. Then $r=2^{t-1}-1$ and $s_p=2^{t-1}$. By Lemma 2.1 (ii), we have $s_p\equiv -1 \ (mod \ p)$. On the other hand, by Lemma 2.1 (i), we get $p|2+2r$ and this leads to $p|1+r$. Then $p|s_p$, a contradiction. 

\textit{Case 7.} Assume that $G\cong C_p^t\rtimes C_2$, where $p$ is an odd prime and $t\geq 1$. Then $s_p=p^t-1$ is even and $s_2>1$ is odd by Lemma 2.1 (ii). Since $\tau_e(G)=\lbrace 1, s_2, s_p\rbrace$ is an arithmetic progression, we obtain $s_p=1+r$ and $s_2=1+2r$. Consequently, $r=p^t-2$. Note that $\pi_e(G)=\lbrace 1, 2, p\rbrace$. Hence $|G|=1+s_2+s_p$ and this leads to $p^t=3$. Then $p=3$, $t=1$ and $G\cong S_3$.

\textit{Case 8.} Let $p, q$ be distinct odd primes. Suppose that $G\cong P\rtimes C_q$, where $P$ is a $p$-group with $exp(P)=p$. Assume that $|P|=p^t$, where $t\geq 1$ is an integer. Then $s_p=p^t-1$, while $s_q=\varphi(q)c_q$. By Lemma 2.1 (iii), we know that $s_p\ne s_q$, so the same-order type of $G$ is $\tau_e(G)=\lbrace 1, s_p, s_q\rbrace$. However, both $s_p$ and $s_q$ are even integers. Consequently, $\tau_e(G)$ is not an arithmetic progression, a contradiction.

\textit{Case 9.} By Theorem 1.2, we assume that $G\cong C_2\times H$, where $H$ is of type (1), (7) or (8).  If $H$ is of type (1) or (8), then $\tau_e(G)=\tau_e(H)$. Consequently, one can replicate the reasoning used in the second paragraph of \textit{Case 1} and the arguments used to solve \textit{Case 8} to arrive at a contradiction. It remains to check what can be said if $H$ is of type (7), i.e. $G\cong C_2\times (C_p^t\rtimes C_2)$, where $p$ is an odd prime and $t\geq 1$. Since $s_p=s_{2p}=p^t-1$ is even, $s_2$ is odd and $\tau_e(G)=\lbrace 1, s_2, s_{p}\rbrace$ is an arithmetic progression, we have $s_p=s_{2p}=1+r$ and $s_2=1+2r$. It follows that $r=p^t-2$. Also, $\pi_e(G)=\lbrace 1, 2, p, 2p \rbrace$, so $|G|=1+s_2+s_p+s_{2p}$. This leads to $r=p^t-1$, a contradiction.    

Conversely, if $G\cong S_3$, we get $\tau_e(G)=\lbrace 1, 2, 3\rbrace$, which is an arithmetic progression of ratio $r=1$. Hence, our proof is complete.
\hfill\rule{1,5mm}{1,5mm}\\

For $k\geq 3$, we denote by $D_{2^k}$ the finite dihedral group with $2^k$ elements. Its well-known structure is given by
$$D_{2^k}=\langle x, y \ | \ x^{2^{k-1}}=y^2=1, x^y=x^{2^{k-1}-1}\rangle.$$
Remark that the result given by Theorem 2.3 holds if $G$ is not a 2-group having more than one cyclic subgroup of order 2. If $G$ is such a group, then $\tau_e(G)$ can be an arithmetic progression. Using GAP \cite{15}, we saw that the first group whose same-order type satisfy this property is $C_8\rtimes C_2^2$ (SmallGroup(32, 43)). Its same-order type is $\lbrace 1, 8, 15\rbrace$. Then, we checked the finite groups of order 64, 128 and 256. For each of these orders, we found 3 examples:
\begin{itemize}
\item[--] order 64: $C_2^2\rtimes D_{16}$ (SmallGroup(64, 128)), $C_2\times(C_8\rtimes C_2^2)$ (SmallGroup(64, 254)), $D_8\circ D_{16}$ (SmallGroup(64, 257)); these groups have the same-order type $\lbrace 1, 16, 31\rbrace$;
\item[--] order 128: $C_2\times (C_2^2\rtimes D_{16})$ (SmallGroup(128, 1728)), $C_2^2\times (C_8\rtimes C_2^2)$ SmallGroup(128, 2310), $C_2\times (D_8\circ D_{16})$ (SmallGroup(128, 2313)); these groups have the same-order type $\lbrace 1, 32, 63\rbrace$;
\item[--] order 256: $C_2^2\times (C_2^2\rtimes D_{16})$ (SmallGroup(256, 53335)), $C_2^3\times (C_8\rtimes C_2^2)$ (SmallGroup(256, 56069)), $C_2^2\times (D_8\circ D_{16})$ (SmallGroup(256, 56072)); these groups have the same-order type $\lbrace 1, 64, 127\rbrace$.
\end{itemize}

Hence, to complete the result given by Theorem 2.3, one would need to solve the following open problem:\\

\textbf{Open problem.} \textit{Classify all finite $2$-groups whose same order types are arithmetic progressions formed of 3 elements.}\\

If we consider the finite groups $G_1\cong C_8\rtimes C_2^2$, $G_2\cong C_2^2\rtimes D_{16}$ and $G_3\cong D_8\circ D_{16}$, then $C_2^t\times G_1$, $C_2^t\times G_2$ and $C_2^t\times G_3$, where $t\geq 0$ is an integer, are for sure included in the answer to the above problem. Their same-order types are arithmetic progressions of the form $\lbrace 1, \frac{|G_i|}{4}, \frac{|G_i|}{2}-1\rbrace$, for $i\in\lbrace 1,2,3\rbrace$.\\

In what follows, we are interested in studying the possibility of finding finite groups whose same-order types are arithmetic progressions formed of 4 elements. Our following main result states there is no finite  nilpotent group $G$ having this property. Firstly, we prove the following preliminary result.\\

\textbf{Lemma 2.4.} \textit{Let $G$ be a finite group of odd order such that $|\tau_e(G)|=4$. Then $\tau_e(G)$ is not an arithmetic progression.}\\

\textbf{Proof.} Let $G$ be a finite group of odd order whose same-order type $\tau_e(G)$ is formed of 4 elements. Then $\tau_e(G)=\lbrace 1, s_l, s_m, s_n\rbrace$, where $l\geq 3, m\geq 3$ and  $n\geq 3$ are distinct odd divisors of $|G|$. Since $\varphi(l), \varphi(m)$ and $\varphi(n)$ are even integers, it follows that $s_l$, $s_m$ and $s_n$ are also even integers. Hence, $\tau_e(G)$ cannot be an arithmetic progression.
\hfill\rule{1,5mm}{1,5mm}\\ 

Next, we consider the case of finite $p$-groups.\\

\textbf{Proposition 2.5.} \textit{Let $G$ be a finite $p$-group such that $|\tau_e(G)|=4$. Then $\tau_e(G)$ is not an arithmetic progression.}\\

\textbf{Proof.} Let $G$ be a finite $p$-group such that $|\tau_e(G)|=4$. If $c_2=1$, by following a similar argument as the one which was written in the first paragraph of \textit{Case 1} in the proof of Theorem 2.3, one obtains $G\cong C_{16}$ or $G\cong Q_{32}$. Then $\tau_e(G)=\lbrace 1, 2, 4, 8\rbrace$ or $\tau_e(G)=\lbrace 1, 4, 8, 16\rbrace$. None of these same-order types is an arithmetic progression. 

Further, we suppose that $c_2>1$.  Assume that $\tau_e(G)$ is an arithmetic progression of ratio $r\geq 1$. Then $\tau_e(G)=\lbrace 1, 1+r, 1+2r, 1+3r\rbrace$. According to Lemma 2.4, it follows that $p=2$, so $G$ is a 2-group of order $2^k$, where $k\geq 4$. Since $s_2=c_2>1$ is odd and $s_4$ is even, the same-order type of $G$ can be written as $\tau_e(G)=\lbrace 1, s_2, s_4, s_{2^a}\rbrace$, where $a\geq 3$ and $s_{2^a}$ is even. Then $s_4=1+r$ and $s_{2^a}=1+3r$ or $s_4=1+3r$ and $s_{2^a}=1+r$. Since $4|s_{2^a}$, we get $4|1+3r$ or $4|1+r$. Then $r=4t+1$ or $r=4t+3$, where $t\geq 0$. If we evaluate $s_4$, we obtain $s_4\equiv 2 \ (mod \ 4)$ in both cases. On the other hand we have $s_4=\varphi(4)c_4$. Consequently, $c_4$ must be an odd integer. According to Theorem 6 of \cite{1}, it follows that $G\cong C_{2^k}, G\cong Q_{2^k}, G\cong D_{2^k}$ or $G\cong S_{2^k}$, where 
$$S_{2^k}=\langle x, y \ | \ x^{2^{k-1}}=y^2=1, x^y=x^{2^{k-2}-1}\rangle$$
is the semidihedral (quasidihedral) group with $2^k$ elements. The first two cases were previously eliminated. If $G\cong D_{2^k}$, then $s_4=2$ and $s_2=2^{k-1}+1$. Since $s_4=1+r$ or $s_4=1+3r$, we obtain $r=1$. Therefore, $s_2=3$ and this leads to $k=2$, a contradiction. If $G\cong S_{2^k}$, then $s_2=2^{k-2}+1$ and $s_4=2^{k-2}+2$. Since $s_2<s_4$, we have $s_4=1+3r$ and $s_{2^a}=1+r$. Consequently, $s_4-s_2=r$ and we deduce that $r=1$. Then $s_{2^a}=2$, but we recall that $4|s_{2^a}$. Hence, we arrive at a contradiction and our proof is finished.
\hfill\rule{1,5mm}{1,5mm}\\

Finally, we extend our previous result to the class of finite nilpotent groups.\\

\textbf{Theorem 2.6.} \textit{Let $G$ be a finite nilpotent group such that $|\tau_e(G)|=4$. Then $\tau_e(G)$ is not an arithmetic progression.}\\

\textbf{Proof.} Let $G$ be a finite nilpotent group such that $|\tau_e(G)|=4$. Assume that $\tau_e(G)$ is an arithmetic progression of ratio $r\geq 1$. Then, according to Lemma 2.4 and Proposition 2.5, $2||G|$ and there is an odd prime $p$ such that $p||G|$. Since $G$ is nilpotent, we have $G\cong H\times K$, where $H$ is the Sylow 2-subgroup of $G$ and $H$ is the direct product of all other Sylow subgroups of $G$. If $s_2=1$, then $H\cong C_2$ and $\tau_e(G)=\tau_e(K)$. In this case, $K$ would be a group of odd order such that its same-order type is an airthmetic progression formed of 4 elements. However, by Lemma 2.4, there are no such groups. Consequently $s_2>1$ and $\tau_e(G)=\lbrace 1, s_2, s_p, s_{2p}\rbrace$. Once again, since $s_2$ is odd, we have $s_2=1+2r$. Also, $s_{2p}=s_2s_p>s_p$, so $s_p=1+r$ and $s_{2p}=1+3r$. Then $1+3r=(1+2r)(1+r)$ and this leads to $r=0$, a contradiction.
\hfill\rule{1,5mm}{1,5mm}\\

Therefore, if there is a finite group $G$ such that $|\tau_e(G)|=4$ and $\tau_e(G)$ is an arithmetic progression, then $G$ is non-nilpotent. We were not able to prove or disprove the existence of such a group. Hence, we add a second open problem.\\

\textbf{Open problem.} \textit{Are there any finite non-nilpotent groups whose same-order types are formed of 4 elements?}\\

According to Lemma 2.4, to provide an answer to the above question it would be sufficient to consider the finite non-nilpotent groups of even order. Using GAP, we checked all such groups $G$ for which $|G|\leq 1000$. None of them satisfies the above arithmetical property. Consequently, according to Lemma 2.1 (iv), if $G$ is a finite group such that $|\tau_e(G)|=4$  and $\tau_e(G)$ is an arithmetic progression of ratio $r$, then $r>3$. Also, it is easy to check that $r$ must be odd. Hence, there are no finite groups whose same-order types are $\lbrace 1, 2, 3, 4\rbrace$ or $\lbrace 1, 4, 7, 10 \rbrace$. Also, in the first section, we recalled that $A_5$ is the only finite non-abelian simple group whose same-order type is formed of 4 elements. Since $\tau_e(A_5)=\lbrace 1, 15, 20, 24\rbrace$ is not an arithmetic progression, it is clear that there are no finite non-abelian simple groups whose same-order types are arithmetic progressions formed of 4 elements.

\vspace*{3ex}
\small

\begin{minipage}[t]{7cm}
Mihai-Silviu Lazorec \\
Faculty of  Mathematics \\
"Al.I. Cuza" University \\
Ia\c si, Romania \\
e-mail: {\tt silviu.lazorec@uaic.ro}
\end{minipage}
\hspace{3cm}
\begin{minipage}[t]{7cm}
Marius T\u arn\u auceanu \\
Faculty of  Mathematics \\
"Al.I. Cuza" University \\
Ia\c si, Romania \\
e-mail: {\tt tarnauc@uaic.ro}
\end{minipage}

\begin{thebibliography}{100}
\bibitem{1} Belshoff, R., Dillstrom, J., Reid, L., \textit{Finite groups with a prescribed number of cyclic subgroups}, Comm. Algebra \textbf{47} (3) (2019), 1043-1056.

\bibitem{2} Berkovich, Y., \textit{Groups of prime power order}, vol. 1, de Gruyter, Berlin (2008).

\bibitem{3} Bianchi, M., Chillag, D., Mauri, A.G.B, Herzog, M., Scoppola, C.M., \textit{Applications of a graph related to conjugacy classes in finite groups}, Arch. Math. \textbf{58} (2) (1992), 126-132.

\bibitem{4} Bianchi, M., Glasby, S.P., Praeger C.E., \textit{Conjugacy class sizes in arithmetic progression}, J. Group Theory \textbf{23} (6) (2020), 1039-1056.

\bibitem{5} Brandl, R., Shi, W., \textit{Finite groups whose element orders are consecutive integers}, J. Algebra \textbf{143} (1991), 388-400.

\bibitem{6} Feng, Y., \textit{Finite groups whose abelian subgroup orders are consecutive integers}, J. Math. Res. Expo. \textbf{18} (1998), 503-506.

\bibitem{7} Frobenius, G., \textit{Verallgeminerung des Sylow'schen Satzes}, Sitzungsberichte
der K\" oniglich Preu\ss ischen Akademie der Wissenschaften zu Berlin (1895), 981-993.

\bibitem{8} Herzog, M., \textit{Counting group elements of order $p$ modulo $p^2$}, Proc. Amer. Math. Soc. \textbf{66} (2) (1977), 247-250.

\bibitem{9} Jafari Taghvasani, L., Zarrin, M., \textit{A characterization of $A_5$ by its same-order type}, Monatsh. Math. \textbf{182} (3) (2017), 731-736.

\bibitem{10} Mazurov, V.D., Khukhro, E.I., \textit{Unsolved Problems in Group Theory. The Kourovka Notebook}, 19th ed., Novosibirsk (2018). 

\bibitem{11} Shen, R., \textit{On groups with given same-order types}, Comm. Algebra \textbf{40} (6) (2012), 2140-2150.

\bibitem{12} Shen, R., Shao, C., Jiang, Q., Shi, W., Mazurov, V.D., \textit{A new characterization of $A_5$}, Monatsh. Math. \textbf{160} (3) (2010), 337-341.

\bibitem{13} Shi, W., \textit{Finite groups whose proper subgroup orders are consecutive integers}, J. Math. Res. Expo. \textbf{14} (1994), 165-166.

\bibitem{14} T\u arn\u auceanu, M., \textit{Arithmetic progressions in finite groups}, arXiv:2003.10060.

\bibitem{15} The GAP Group, \textit{GAP -- Groups, Algorithms, and Programming, Version 4.9.3}, 2018, https://www.gap-system.org.
\end{thebibliography}
\end{document}